\documentclass[12pt]{article}
\usepackage{amsmath, amsthm, amsfonts}
\usepackage{amssymb}
\usepackage{amscd}
\usepackage{verbatim}

\begin{document}

\newcommand{\merhav}{{\mathcal D}^{1,2}}
\newcommand{\be}{\begin{equation}}
\newcommand{\ee}{\end{equation}}
%%%%%%%%%%%%%%%
\newcommand{\bea}{\begin{eqnarray}}
\newcommand{\eea}{\end{eqnarray}}
\newcommand{\bean}{\begin{eqnarray*}}
\newcommand{\eean}{\end{eqnarray*}}
\newcommand{\thkl}{\rule[-.5mm]{.3mm}{3mm}}
%%%%%%%%%%%%%%%%%%%%%%%%%%%
\newcommand{\cw}{\stackrel{\rightharpoonup}{\rightharpoonup}}
\newcommand{\id}{\operatorname{id}}
\newcommand{\supp}{\operatorname{supp}}
\newcommand{\wlim}{\mbox{ w-lim }}
\newcommand{\mymu}{{x_N^{-p_*}}}
\newcommand{\R}{{\mathbb R}}
\newcommand{\N}{{\mathbb N}}
\newcommand{\Z}{{\mathbb Z}}
\newcommand{\Q}{{\mathbb Q}}
\newtheorem{theorem}{Theorem}[section]
\newtheorem{corollary}[theorem]{Corollary}
\newtheorem{lemma}[theorem]{Lemma}
\newtheorem{definition}[theorem]{Definition}
\newtheorem{remark}[theorem]{Remark}
\newtheorem{proposition}[theorem]{Proposition}
%%%%%%%%%%%%%%%%%%
\newtheorem{conjecture}[theorem]{Conjecture}
\newtheorem{question}[theorem]{Question}
\newtheorem{example}[theorem]{Example}
\newtheorem{Thm}[theorem]{Theorem}
\newtheorem{Lem}[theorem]{Lemma}
\newtheorem{Pro}[theorem]{Proposition}
\newtheorem{Def}[theorem]{Definition}
\newtheorem{Exa}[theorem]{Example}
\newtheorem{Exs}[theorem]{Examples}
\newtheorem{Rems}[theorem]{Remarks}
\newtheorem{Rem}[theorem]{Remark}
\newtheorem{Cor}[theorem]{Corollary}
\newtheorem{Conj}[theorem]{Conjecture}
\newtheorem{Prob}[theorem]{Problem}
\newtheorem{Ques}[theorem]{Question}
%%%%%%%%%%%%%%%%%%
%\newenvironment{proof}{{\bf Proof.}}{\hfill $\bowtie$\vskip4mm}

\renewcommand{\theequation}{\thesection.\arabic{equation}}
\catcode`@=11
\@addtoreset{equation}{section}
\catcode`@=12

%\begin{titlepage}

\title{Existence of minimizers for Schr\"{o}dinger operators under domain
perturbations with application to Hardy's inequality}
\author{Yehuda Pinchover\\
 {\small Department of Mathematics}\\ {\small  Technion - Israel Institute of Technology}\\
 {\small Haifa 32000, Israel}\\
{\small pincho@techunix.technion.ac.il}\\\and Kyril Tintarev
\\{\small Department of Mathematics}\\{\small Uppsala University}\\
{\small SE-751 06 Uppsala, Sweden}\\{\small
kyril.tintarev@math.uu.se}}
%\date{}
\maketitle
\newcommand{\dnorm}[1]{\thkl #1 \thkl\,}

\begin{abstract} The paper studies the existence of minimizers for
Rayleigh quotients $\mu_{\Omega}=\inf\frac{\int_\Omega|\nabla
u|^2}{\int_\Omega V{|u|^2}}\,$, where $\Omega$ is a domain in
$\mathbb{R}^N$, and $V$ is a nonzero nonnegative function that may
have singularities on $\partial\Omega$. As a model for our results
one can take $\Omega$ to be a Lipschitz cone and $V$ to be the
Hardy potential $V(x)=\frac{1}{|x|^2}\,$.
\\[2mm]
\end{abstract}

%\end{titlepage}
%%%%%%%%%%%%%%%%%%%%%%%%%%%%%%%%%%%%%%%%%%%%%%%%%
\section{Introduction}
Let $X$ be a domain in $\mathbb{R}^N$, and let $V\in L^p_{loc}(X)$
be a nonzero nonnegative function, where $p>\frac{N}{2}$. Let
$\merhav(X)$ be the completion of $C_0^\infty(X)$ with respect to
the norm $\|u\|^2=\int_X|\nabla u|^2$. For an open set
$\Omega\subset X$, we will consider the subspace
$\merhav(\Omega)\subset \merhav(X)$, which is by definition, the
closure in $\merhav(X)$ of $C_0^\infty(\Omega)$. We denote
$B\Subset X$, if $B\subset X$,  and $\overline{B}$ is compact in
$X$.

Let $\Omega\subset X$. We study the existence of a minimizer for
the Rayleigh quotient \be \label{mu} \mu_{\Omega}=\inf_{u\in
\merhav(\Omega),\,Vu\not\equiv 0} \frac{\int_\Omega|\nabla
u|^2}{\int_\Omega V{|u|^2}}\,,\ee under the assumption that
 \be \label{Hardy} \mu_X>0. \ee
Condition (\ref{Hardy}) is satisfied, for example, when
$X\!=\!\R^N\!\setminus\!\{0\}$, $V(x)\!=\!\frac{1}{|x|^2}\,$, and
$N\geq 3$, which corresponds to the well-known Hardy inequality,
with $\mu_X=\frac{(N-2)^2}{4}$. Existence of a minimizer  in
problems with a singular potential has been studied by many
authors with attention to `small' perturbations of the potential
$V$ (see, \cite{BM,BMS,CH,MMP,MS,S,T} and the references therein).
Typically in such cases, if there is a `spectral gap', then a
minimizer exists. This situation is sometimes called the `gap
phenomenon'. The present paper studies the existence of a
minimizer in the case of compact domain perturbations under the
situation of a positive `spectral gap'. Domain perturbations in
the context of variational inequalities and the Dirichlet problem
were studied in \cite{D,Mosco} and the references therein.

Let $P$ be a second order elliptic operator which is defined on a
domain $\Omega$, and denote by $\mathcal{C}_P(\Omega)$ the cone of
all positive solutions of the equation $Pu=0$ in $\Omega$. For
$P_\mu:=-\Delta-\mu V$, we simply write
$\mathcal{C}_\mu(\Omega):=\mathcal{C}_{P_\mu}(\Omega)$. Let
$K\Subset \Omega$. Recall \cite{P,Pinsk} that $u\in
\mathcal{C}_P(\Omega\setminus K)$ is said to be a {\em positive
solution of the operator $P$ of minimal growth in a neighborhood
of infinity in} $\Omega$, if for any $K\Subset K_1 \Subset \Omega$
and any $v\in C(\overline{\Omega\setminus K_1})\cap
\mathcal{C}_P(\Omega\setminus K_1)$, the inequality $u\le v$ on
$\partial K_1$ implies that $u\le v$ in $\Omega\setminus K_1$. A
positive solution $u\in \mathcal{C}_P(\Omega)$ which has minimal
growth in a neighborhood of infinity in $\Omega$ is called a {\em
ground state} of $P$ in $\Omega$.

The operator $P$ is said to be {\em critical} in $\Omega$, if $P$
admits a ground state in $\Omega$. The operator $P$ is called {\em
subcritical} in $\Omega$, if $\mathcal{C}_P(\Omega)\neq
\emptyset$, but $P$ is not critical in $\Omega$. If
$\mathcal{C}_P(\Omega)= \emptyset$, then $P$ is {\em
supercritical} in $\Omega$.

Suppose that $P$ is critical in $\Omega\varsubsetneq X$. Then $P$
is subcritical in any domain $\Omega_1$ such that
$\Omega_1\varsubsetneq \Omega$, and supercritical in any domain
$\Omega_2$ such that $\Omega\varsubsetneq \Omega_2\subset X$.
Furthermore, for any nonzero nonnegative function $W$ the operator
$P+W$ is subcritical and $P-W$ is supercritical in $\Omega$.
Moreover, if $P$ is critical in $\Omega$, then $\dim
\mathcal{C}_P(\Omega)=1$ (see e.g. \cite{Pinsk}).

If $P$ is subcritical in $\Omega$, then $P$ admits a positive
minimal Green function $G^\Omega_P(x,y)$ in $\Omega$. Moreover,
for each $y\in \Omega$, the function $G^\Omega_P(\cdot,y)$ is a
positive solution of the equation $Pu=0$ in $\Omega\setminus\{y\}$
that has minimal growth in a neighborhood of infinity in $\Omega$
(see \cite{Pinsk}).

Consider now the case that $P_\mu=-\Delta-\mu V$, where $V$ is a
nonzero nonnegative function and $\mu\in \mathbb{R}$. Then $P_\mu$
is subcritical in $\Omega$ for all $\mu<\mu_\Omega$, supercritical
in $\Omega$ for all $\mu>\mu_\Omega$, and $P_{\mu_\Omega}$ is
either critical or subcritical, where $\mu_\Omega$ is defined by
(\ref{mu}).

In many papers the term {\em ground state} refers only to
minimizer solutions of (\ref{mu}). It turns out that such a
minimizer solution  is also the ground state of the operator
$-\Delta-\mu_\Omega V$ in the sense introduced above. For
Schr\"{o}dinger operators this fact was proved in \cite{Mu0} (see
Theorem 2.7 therein, and the remark below its proof). The
following lemma applies also to the general symmetric case, and
its proof applied even to nonsymmetric cases. An alternative proof
that was suggested to us by M.~Murata (after the first draft of
the present paper has been completed) uses the heat kernel.

\begin{lemma}\label{lemminimizer} Suppose that $V>0$ and (\ref{mu}) admits a
minimizer, then the operator $-\Delta-\mu_\Omega V$ is critical in
$\Omega$, and a minimizer is a ground state.
\end{lemma}

Our first main result reads as follows.
\begin{theorem} \label{thm1subc} Suppose
that $\Omega\subset X$ is a domain satisfying
$0<\mu_X<\mu_\Omega$. Then there exists an open set $B\Subset X$
such that $\Omega\cup B$ is connected and \be \label{unequal mu}
\mu_{\Omega\cup B}<\mu_\Omega.\ee Moreover, for any such set $B$
the infimum value $\mu_{\Omega\cup B}$ for problem (\ref{mu}) is
uniquely attained.
\end{theorem}
\begin{corollary}
Suppose that $B$ satisfies the conditions of
Theorem~\ref{thm1subc}. Then for every open set $B^\prime$ such
that $B\subset B^\prime\Subset X$ and $\Omega\cup B^\prime$ is
connected, the infimum value $\mu_{\Omega\cup B^\prime}$ is
attained.
\end{corollary}
\begin{proof} The inclusion  $\merhav(\Omega\cup B)\subset
\merhav(\Omega\cup B^\prime)$ implies $\mu_{\Omega\cup
B^\prime}\le \mu_{\Omega\cup B}<\mu_\Omega$, and hence
Theorem~\ref{thm1subc} applies.
\end{proof}
In the critical case we have the following stronger statement.
\begin{theorem} \label{thm1cr} Suppose
that $\Omega\subset X$ is a domain satisfying
$0<\mu_X<\mu_\Omega$, and assume that the operator
$P=-\Delta-\mu_\Omega V$ is critical in $\Omega$.

\vskip 2mm

 \noindent Then for any open set $B\Subset
X$ such that $\Omega\cup B$ is connected and $\Omega\neq\Omega\cup
B$, the inequality (\ref{unequal mu}) is satisfied, and the
infimum value $\mu_{\Omega\cup B}$ for problem (\ref{mu}) is
uniquely attained.
\end{theorem}
If $\Omega\Subset X$, then it is well known that $\mu_\Omega$ in
(\ref{mu}) is attained since $\int_\Omega V|u|^2$ is weakly
continuous. For a noncompact domain $\Omega$,  or a potential $V$
that blows up near $\partial \Omega\cap \partial X$ the minimizer
may not exist, as the following example demonstrates.
\begin{example}\label{excone}{\em
Consider a Lipschitz (connected) cone $\mathrm{C}\subset
\R^N\setminus\{0\}$, $N\geq 2$, with the vertex at $0$. Let
$V(x)=\frac{1}{|x|^2}$, and $\mu\in \mathbb{R}$. Denote by
$\mathcal{C}^0_\mu(\mathrm{C})$ the cone of all positive solutions
of the equation \be \label{eigenfunctions} P_\mu u: =-\Delta u-
\mu\frac{u}{|x|^2}=0 \ee in $\mathrm{C}$ that vanish on $\partial
\mathrm{C}\setminus \{0\}$. By \cite{P}, the dimension of
$\mathcal{C}^0_\mu(\mathrm{C})$ is at most 2. Actually, using
separation of variables and \cite{P}, one can compute the
solutions in $\mathcal{C}^0_\mu(\mathrm{C})$ explicitly.

Let $D\subset S_1^{N-1}$ be the Lipschitz domain so that
$$\mathrm{C}=\{(r,\omega)\mid \,r\in(0,\infty), \omega\in D\}.$$
Denote by $\Delta_r$ and $\Delta_S$ the radial and the spherical
Laplacian, respectively. Let $\lambda_D$ and $v_D(\omega)$ be the
Dirichlet principal eigenvalue and eigenfunction of $-\Delta_S$ on
$D$. So,
$$-\Delta_Sv_D=\lambda_D v_D \quad \mbox{ on } D, \qquad  v_D\!\left|_{\partial D}\right. =0 .$$
Then any positive solution in $\mathcal{C}^0_\mu(\mathrm{C})$ is
of the form
$$u_{\mu,D}(r)v_D(\omega)\qquad   r\in(0,\infty),\, \omega\in
S_1^{N-1},$$ where $u_{\mu,D}$ is a {\em global} positive solution
of the Euler equidimensional equation
$$-\Delta_r u-\frac{\mu-\lambda_D}{r^2}u=-u''-\frac{(N-1)}{r}u'-\frac{\mu-\lambda_D}{r^2}u=0\qquad 0<r<\infty.$$
It follows that $\mu$ should satisfy $\mu\leq
\frac{(N-2)^2}{4}+\lambda_D$, and
$u_{\mu,D}(r)=ar^{\alpha_+}+br^{\alpha_-}$, where
$$\alpha_\pm=\alpha_\pm(\mu,D)=\frac{-(N-2)\pm\sqrt{(N-2)^2-4(\mu-\lambda_D)}}{2}\,,$$
and $a,b\geq0$. In particular, \be\label{muc}
\mu_\mathrm{C}=\frac{(N-2)^2}{4}+\lambda_D,\ee
 and the
corresponding unique positive solution in
$\mathcal{C}^0_{\mu_\mathrm{C}}\!(\mathrm{C})$ equals
$r^{-\frac{(N\!-\!2)}{2}}v_D(\omega)$, which clearly does not
belong to $\merhav(\mathrm{C})$. It is well known that if a
minimizer of the variational problem exists, then it belongs to
$\mathcal{C}^0_{\mu_\mathrm{C}}(\mathrm{C})$. Therefore,
$\mu_\mathrm{C}$ is not attained for any Lipschitz cone
$\mathrm{C}\subset \mathbb{R}^N$. On the other hand, noting that
the solution $ r^{-\frac{(N-2)}{2}}\log r\,v_D(\omega)$ is a
positive solution of the equation $P_{\mu_\mathrm{C}} u=0$ near
$\zeta=0$ and $\zeta=\infty$ which grows there faster than
$r^{-\frac{(N-2)}{2}}v_D(\omega)$, and using \cite{P}, it follows
that $r^{-\frac{(N-2)}{2}}v_D(\omega)$ is a ground state of the
{\em critical} operator $P_{\mu_\mathrm{C}}$ in $\mathrm{C}$.

Now, for $N\geq 3$ take $X:=\mathbb{R}^N\setminus \{0\}$, and note
that  $\mu_X=\frac{(N-2)^2}{4}>0$, so, (\ref{Hardy}) is satisfied.
For $N=2$ take a Lipschitz cone $X$ with a vertex at the origin
such that $\overline{\mathrm{C}}\setminus \{0\}\subset
X\varsubsetneq \R^N\setminus\{0\}$. So, (\ref{Hardy}) is satisfied
also in the two dimensional case.

Consequently, Theorem~\ref{thm1cr} implies that  for any open set
$B\Subset X$ such that $\mathrm{C}\cup B$ is connected, and
$\mathrm{C}\varsubsetneq \mathrm{C}\cup B$, the infimum value
$\mu_{\mathrm{C}\cup B}$ is uniquely attained. By \cite{P}, it
follows that the corresponding minimizer behaves near
$\zeta=\infty$ and near $\zeta=0$ like
$r^{\alpha_-(\mu_{\mathrm{C}\cup B},D)}v_D(\omega)$ and
$r^{\alpha_+(\mu_{\mathrm{C}\cup B},D)}v_D(\omega)$, respectively.

On the other hand, if $B$ is replaced by a larger set that is not
relatively compact in $X$, then a minimizer may not exist. Take
for example two connected Lipschitz cones $\mathrm{C}$ and
$\mathrm{C}_1$, such that $\mathrm{C}\varsubsetneq
\mathrm{C}_1\subset X$. Notice that one has
$\lambda_{D_1}<\lambda_{D}$, and by (\ref{muc}),
$\mu_{\mathrm{C}_1}<\mu_{\mathrm{C}}$. Hence, for $B=\mathrm{C}_1$
we have $\mathrm{C}\cup B=\mathrm{C}_1$, and consequently, the
infimum $\mu_{\mathrm{C}\cup B}$ is not attained.
 }\end{example}

Next, we discuss the subcritical case, where adding a compact set
that is too small, also implies the  non-existence of a minimizer:

\begin{theorem}\label{thm2}
Let $\Omega\varsubsetneq X$ be a domain with a Lipschitz boundary,
and let $V\in C^\alpha_{loc}(X)$ be a positive function, where
$0<\alpha\leq 1$. Assume that the operator $P:=-\Delta-\mu_\Omega
V$ is subcritical in $\Omega$, and (\ref{Hardy}) is satisfied.

Let $B_j\Subset X$ be a decreasing sequence of smooth domains,
such that $\Omega_j:=B_j\cup \Omega$ are connected for all $j\geq
1$, ${\rm int}\left(\cap_j\Omega_j\right)=\Omega$, and $B_1
\cap\partial \Omega$ is contained in a Lipschitz portion $\Gamma
\Subset\partial \Omega$. Then there exists $j_0>0$ such that for
all $j\ge j_0$, $\mu_{\Omega_j}$ is not attained. Moreover,
$-\Delta-\mu_\Omega V$ is subcritical in $\Omega_j$ for all $j\ge
j_0$.
\end{theorem}

In particular, we have
\begin{corollary}\label{thm2old}
Let $\mathrm{C}\varsubsetneq X$ be a Lipschitz cone with vertex at
$0$, where $X=\R^N\setminus\{0\}$ if $N\geq 3$, and
$X\varsubsetneq \R^N\setminus\{0\}$ is a Lipschitz cone with a
vertex at the origin such that $\overline{\mathrm{C}}\setminus
\{0\}\subset X$, if $N=2$. Let $W\in C^\alpha_{loc}(X)$,
$0<\alpha\leq 1$ be a nonzero nonnegative function with a compact
support in $\mathrm{C}$, and set $V(x)=\frac{1}{|x|^2}-W(x)$. Let
$B_j\Subset X$ be a decreasing sequence of smooth domains, such
that $\mathrm{C}_j:=B_j\cup \mathrm{C}$ are connected for all
$j\geq 1$, and ${\rm
int}\left(\cap_j\mathrm{C}_j\right)=\mathrm{C}$. Then there exists
$j_0>0$ such that for all $j\ge j_0$, $\mu_{\mathrm{C}_j}$ is not
attained, and $-\Delta-\mu_\mathrm{C} V$ is subcritical in
$\mathrm{C}_j$ for all $j\ge j_0$.
\end{corollary}
%%%%%%%%%%%%%%%%%%%%%%%%%%%%%%%%%%%%%%%%%%%%%%%%%%%%%%%%%%%
\section{Existence of minimizers under compact domain perturbations}
In this section we prove Theorem~\ref{thm1subc} and
Theorem~\ref{thm1cr}. Throughout the section we assume that
$\mu_X<\mu_\Omega$.
\begin{lemma}
\label{B_epsilon} For any $\varepsilon>0$ there exists an open
bounded set $B_\varepsilon\Subset X$, such that
$\mu_{B_\varepsilon}\le\mu_X+\varepsilon$.
\end{lemma}

\begin{proof}
Since $C_0^\infty(X)$ is dense in $\merhav(X)$, there exists a
minimizing sequence $u_k\in C_0^\infty(X)$, such that $\int_X
V|u_k|^2=1$ and $\|u_k\|^2\le\mu_X+k^{-1}$. Fix
$k_\varepsilon>\varepsilon^{-1}$, and choose an open bounded set
$B_\varepsilon$ so that $\supp u_{k_{\varepsilon}}\subset
B_\varepsilon\Subset X$. Then $u_{k_\varepsilon}\in
\merhav(B_\varepsilon)$ and
$\mu_{B_\varepsilon}\le\|u_{k_\varepsilon}\|^2\le\mu_{X}+k^{-1}_\varepsilon<\mu_X+\varepsilon$.
\end{proof}
 Let $0<\varepsilon<\mu_\Omega-\mu_X$. Since
$\merhav(B_\varepsilon)\subset\merhav(\Omega\cup B_\varepsilon)$,
we have  $$\mu_{\Omega\cup
B_\varepsilon}\le\mu_{B_\varepsilon}<\mu_X+\varepsilon<\mu_\Omega\,.$$
Recall that if the operator $P=-\Delta-\mu_\Omega V$ is critical
in $\Omega$, and $\Omega\varsubsetneq\Omega_1$, then
$\mu_{\Omega_1}<\mu_{\Omega}$. Consequently, the assertions of
theorems \ref{thm1subc} and \ref{thm1cr} follow from the following
statement.

\begin{lemma} If $B\Subset X$ is an open set, and
$\mu_{\Omega\cup B}<\mu_\Omega$, then $\mu_{\Omega\cup B}$ is
attained and every minimizing sequence for $\mu_{\Omega\cup B}$ is
convergent.
\end{lemma}

\begin{proof} Let $\{u_k\}$ be a minimizing sequence for $\mu_{\Omega\cup
B}$. So, we may assume that $\int_{\Omega\cup B}V|u_k|^2=1$ and
$\|u_k\|^2\to\mu_{\Omega\cup B}$. Consider a weakly convergent in
$\merhav(\Omega\cup B)$  subsequence of $\{u_k\}$, which we
relabel as $\{u_k\}$. Let $w:=\wlim u_k$, and denote
$v_k:=u_k-w\rightharpoonup 0$. Since $(v_k,w)\to 0$, we have
\be\|u_k\|^2=\|v_k+w\|^2=\|v_k\|^2+\|w\|^2+2{\rm Re}\,(v_k,w) =
\|v_k\|^2+\|w\|^2+o(1),\ee so that \be \label{energy}
\|v_k\|^2+\|w\|^2=\mu_{\Omega\cup B}+o(1).\ee

Note that $\int_{\Omega\cup B}Vv_kw\to 0$, since (\ref{Hardy}) and
Cauchy-Schwartz inequality imply that $u\mapsto\int_{\Omega\cup
B}Vuw$ is a continuous functional on $\merhav(\Omega\cup B)$.
Thus, by repeating the derivation of (\ref{energy}) for the
seminorm $\sqrt{\int V|u|^2}$, we have,  \be \label{mass}
\int_{\Omega\cup B} V|v_k|^2+\int_{\Omega\cup B} V|w|^2 =1+o(1).
\ee Let $t=\int_{\Omega\cup B} Vw^2 $. Once we show that
 \be\label{crucial}\|v_k\|^2\geq\mu_\Omega\int_{\Omega\cup B} V|v_k|^2
+o(1),\ee we will have from (\ref{energy}) and (\ref{mass}) that
$(1-t)\mu_\Omega+t\mu_{\Omega\cup B}\leq \mu_{\Omega\cup B}$.
Since $\mu_\Omega>\mu_{\Omega\cup B}$, this can hold only if
$t=1$. By (\ref{energy}),  $\mu_{\Omega\cup B}\ge\|w\|^2$ and
since $\int_{\Omega\cup B} V|w|^2=1$ we see that $w$ is a
minimizer. Moreover, since $\|u_k\|\to\|w\|$, $u_k\to w$ in
$\merhav$.

Let us verify (\ref{crucial}). Let $\chi\in C_0^\infty(X;[0,1])$
be equal 1 on $\overline{B}$. Then,  by the compactness of the
Sobolev imbedding on bounded smooth sets, we have
$\int_{\supp\chi}V|v_k|^2\to 0$, and

\be \label{mass1} \int_{\Omega\cup B} V|v_k|^2 = \int_{\Omega\cup
B} V[(1-\chi)^2 |v_k|^2+\chi(2-\chi)|v_k|^2] = \int_{\Omega}
V|(1-\chi)v_k|^2+o(1).\ee Observe that
\begin{eqnarray}
 \nonumber
\int_{\Omega\cup B}|\nabla v_k|^2-\int_{\Omega\cup B}|\nabla\left((1-\chi) v_k\right)|^2=\\
\nonumber -\int_{\Omega\cup
B}\left(|\nabla\chi|^2|v_k|^2-2(1-\chi) v_k\nabla\chi\cdot\nabla
v_k\right)
+\\
\label{en2} \int_{\Omega\cup B}\chi(2-\chi)|\nabla v_k|^2.
\end{eqnarray}
By the compactness of Sobolev imbedding on relatively compact
smooth sets, we have \be \label{en3} \int_{\Omega\cup
B}|\nabla\chi|^2|v_k|^2\le C\int_{\supp\chi}|v_k|^2=o(1),\ee
 and \be \label{en4} \left|\int_{\Omega\cup
B}(1-\chi) v_k\nabla\chi\cdot\nabla v_k\right|\le C
\left(\int_{\Omega\cup B}|\nabla
v_k|^2\right)^\frac12\left(\int_{\supp\chi}|v_k|^2\right)^\frac12
=o(1).\ee
 Combining (\ref{en2}),(\ref{en3}) and (\ref{en4}), we have
\be\label{e58}\|v_k\|^2\ge \int_{\Omega\cup
B}|\nabla\left((1-\chi) v_k\right)|^2+o(1). \ee

 \noindent{\bf Claim:} For any $\psi\in \merhav(\Omega\cup B)$, we have
 $(1-\chi)\psi\in\merhav(\Omega)$.  Let $\{\psi_l\}_{l=1}^\infty \subset
C_0^\infty (\Omega\cup B)$ be a sequence such that $\psi_l\to
\psi$ in $\merhav(\Omega\cup  B)$. Since $(1-\chi)\psi_l \in
\merhav(\Omega)$, it is enough to show that $(1-\chi)\psi_l \to
(1-\chi)\psi$ in $\merhav(\Omega)$.

 Indeed
\begin{eqnarray*}
\int_{\Omega}|\nabla \left((1-\chi)(\psi_l -\psi)\right)|^2\leq\\
 2\int_{\Omega}|\nabla(1-\chi)|^2|\psi_l-\psi|^2+
2\int_{\Omega} |(1-\chi)|^2|\nabla(\psi_l-\psi)|^2\leq\\
 2\int_{\Omega}|\nabla(1-\chi)|^2|\psi_l-\psi|^2+
2\int_{\Omega\cup B} \nabla(\psi_l-\psi)|^2 \to 0,
\end{eqnarray*}
where we used the compactness of Sobolev imbedding on relatively
compact smooth sets.

  By the Claim
$(1-\chi)v_k\in\merhav(\Omega)$, therefore, (\ref{e58}) and the
definition of $\mu_\Omega$  imply
 \be\|v_k\|^2\ge
\mu_\Omega\int_{\Omega}V|(1-\chi) v_k|^2+o(1).\ee
 Substituting (\ref{mass1}) into the last inequality,
 we obtain (\ref{crucial}), which proves
the lemma.
\end{proof}
%%%%%%%%%%%%%%%%%%%%%%%%%%%%%%%%%%%%%%%%%%%%%%%%%%%%%%%%%%%%%%
\section{Proof of Lemma \ref{lemminimizer}}
Throughout this section,  $\Omega$ denotes a domain  in
$\mathbb{R}^N$, $N\geq 2$, and $V>0$. We start with a brief
discussion of some spectral properties of the operator
$P=-V(x)\Delta$ in $\Omega$.

 First, we turn $\Omega$ into a Riemannian manifold $M$
 equipped with the metric $ds^2= (V(x))^{-1} \sum_{i=1}^N dx_i^2$.
We put $\tilde{L}_2(M)=L_2(\Omega;V)$ equipped with the norm
$\dnorm{u}_2=(\int_\Omega |u|^2Vdx)^{\frac{1}{2}}$,  and
$$\tilde{H}^{1}(M)= \{u\in W_{1,2}^{loc}(\Omega)\,:\,
\dnorm{u}_{1,2}:=\left(\dnorm{u}_2^2+\|\nabla
u\|_{L_2(\Omega)}^2\right)^{1/2}<\infty\}.$$
 The closure of $C^1_0(\Omega)$ under this norm
will be denoted by $\tilde{H}^{1}_0(M)$.

 Let $\tilde{P}$ be the Friedrichs extension  of the
operator $P$ considered as a symmetric operator in
$\tilde{L}_2(M)$ with domain $C^1_0(\Omega)$ (see \cite{A0}).
\begin{remark}\label{remcomplete}{\em
If $M$ is a complete Riemannian manifold, then the operator
$\tilde{P}$ is the unique selfadjoint realization of $P$ in
$\tilde{L}_2(M)$. In this case, $\tilde{P}$ coincides with the
Dirichlet realization of $P$ with domain of definition given by
$$D(\tilde{P})=\{u\; | \; u \in \tilde{L}_2(M) \cap H^1_{loc}(M), Pu \in
\tilde{L}_2(M)\}.$$
 }\end{remark}

\vskip 3mm

 We denote by $\sigma(\tilde{P})$, $\sigma_{point}(\tilde{P})$, the
spectrum  and point spectrum  of $\tilde{P}$, respectively.

It is well known that
$$
\lambda_0:=\inf \sigma(\tilde{P})=\mu_\Omega= \inf_{u\in
\tilde{H}^1_0(M)} \frac{\int_{\Omega}|\nabla u|^2 dx}{\int_
{\Omega} |u|^2V dx}\;,$$ and  \bean\lambda_0 &=& \sup\{\lambda \in
\mathbb{R} \; :\;
\mathcal{C}_{P-\lambda}(\Omega)\neq \emptyset \}\\
&=&\sup\{\lambda \in \mathbb{R} \; : \; \exists u \in
H^{1}_{loc}(\Omega), u>0, (P-\lambda)u \geq 0 \mbox{ in } \Omega
\}, \eean
 and the supremum $\lambda_0$ is achieved.

 \begin{comment}Further,
\be \label{varp1} \lambda_{\infty}:= \inf \sigma_{ess}(\tilde{P})=
\sup_{K\Subset \Omega}\lambda_0(M\setminus K), \ee where
$\lambda_0(M\setminus K)$ is defined in the same way as
$\lambda_0$ with $M$ replaced by $M\setminus K$, and \bean
 \lambda_{\infty}&=& \sup\{\lambda \in \mathbb{R}
\;:\; \exists K\subset \subset \Omega \mbox{ s.t.  }
{\cal C}_{P-\lambda}(\Omega \setminus K)\neq \emptyset \}\\
&=&\sup\{\lambda \in \mathbb{R} \; : \;\exists K\subset \subset
\Omega \mbox{ and }
u \in H^{1}_{loc}(\Omega\setminus K), \mbox{ s.t. } \\
& & \hspace{2.2cm} u>0 \;\mbox{ and }\; (P-\lambda)u \geq 0
\;\mbox{ in } \;\Omega \setminus K\}, \eean (see \cite{Ag}).
Clearly, $\lambda_{\infty} \geq \lambda_0$.
\end{comment}

If the infimum in (\ref{mu}) is achieved, then it possesses a
positive minimizer. Since every minimizer is a solution of the
equation $(P-\lambda_0)u=0$ in $\Omega$, it follows that problem
(\ref{mu}) possesses a minimizer $\varphi$ if and only if
$\lambda_0=\mu_\Omega\in \sigma_{point}(\tilde{P})$ and
$\varphi\in \mathcal{C}_{P-\lambda_0}(\Omega) \cap
\tilde{L}_2(\Omega)$.

\begin{proof}[Proof of Lemma \ref{lemminimizer}]
 By the Birman-Schwinger principle, $\lambda_0\in
\sigma_{point}(\tilde{P})$ if and only if there exists $\varphi\in
\tilde{L}_2(M)$ such that for every $0\leq \lambda<\lambda_0$ we
have in the $L_2$ sense
\begin{equation}\label{BSCH}
(\lambda_0-\lambda)\int_\Omega V(x)^{1/2}G_{-\Delta-\lambda
V}^\Omega(x,y)V(y)^{1/2}\left(V(y)^{1/2}\varphi(y)\right)\,dy=
V(x)^{1/2}\varphi(x).
\end{equation}
Moreover, by the continuity of the minimizer $\varphi$ and the
positivity of $V$, (\ref{BSCH}) holds
 true if and only if
\begin{equation}\label{BSCH4}
\int_\Omega G_{-\Delta-\lambda V}^\Omega(x,y)V(y)\varphi(y)\,dy=
\frac{\varphi(x)}{\lambda_0-\lambda}
\end{equation}
for all $x\in \Omega$.

 Fix $x_0\in \Omega$.
Assume that $-\Delta-\lambda_0 V$ is subcritical in $\Omega$.
Since $\varphi\in \mathcal{C}_{P-\lambda_0}(\Omega)$ and
$G_{-\Delta-\lambda_0 V}^\Omega(x_0,x)$ is a positive solution of
the operator $-\Delta-\lambda_0 V$ of minimal growth in a
neighborhood of infinity in $\Omega$, it follows that there exists
$C_\varepsilon>0$ such that $G_{-\Delta-\lambda_0
V}^\Omega(x_0,x)\leq C_\varepsilon\varphi(x)$ for all $x\in
\Omega\setminus B(x_0,\varepsilon)$. In particular,
$$G_{-\Delta-\lambda_0
V}^\Omega(x_0,y)V(y)\varphi(y)\in L_1(\Omega).$$ By the Lebesgue
monotone convergence theorem \bean\label{BSCH7} \int_\Omega
G_{-\Delta\!-\!\lambda_0
V}^\Omega(x_0,y)V(y)\varphi(y)\,dy=\\[2mm]\lim_{\lambda\nearrow
\lambda_0}\int_\Omega G_{-\Delta-\lambda
V}^\Omega(x_0,y)V(y)\varphi(y)\,dy= \lim_{\lambda\nearrow
\lambda_0}\!\frac{\varphi(x_0)}{\lambda_0-\lambda}=\infty, \eean
which is a contradiction. Therefore $-\Delta-\lambda_0 V$ is
critical, and $\varphi$ is a ground state of the operator
$-\Delta-\lambda_0 V$ in $\Omega$.
\end{proof}
\section{Nonexistence of minimizers under small\\ compact domain
perturbations
 %\\in the subcritical case
 } In this section we prove
Theorem~\ref{thm2} and give a direct proof of
Corollary~\ref{thm2old}.
\begin{proof} ({\bf proof of Theorem \ref{thm2}})
Consider the domain $\Omega\varsubsetneq X$, and   let $B_j\Subset
X$ be the given decreasing sequence. Consider the  Lipschitz
portion $\Gamma\subset\partial \Omega$ such that $B_1 \cap\partial
\Omega$ is contained in $\Gamma$.

Let $\Gamma_\varepsilon$ denote the set
$$\Gamma_\varepsilon:=\{x\in \Omega\mid
\mbox{dist\,}(x,\Gamma)=\varepsilon\},$$ where $\varepsilon>0$ is
sufficiently small. Finally, fix $x_0\in \Omega$ such that
$\mbox{dist\,}(x_0,\Gamma)=\varepsilon/2$.

Suppose that that $\mu_j:=\mu_{\Omega_j}$ is attained for all
$j\ge 1$, and let $u_j\in \mathcal{C}_{\mu_j}(\Omega_j) \cap
\merhav(\Omega_j) $ be the corresponding minimizer such that
$u_j(x_0)=1$.  By Lemma \ref{lemminimizer}, $u_j$ is  the
normalized ground state of the (critical) operator
$P_j:=-\Delta-\mu_j V$ in $\Omega_j$.

Clearly,  $\mu_j\leq \mu_\Omega$. Therefore, $P_j$ is subcritical
in $\Omega$, and denote by $G_{P_j}^\Omega(x,x_0)$  the
corresponding positive minimal Green function.

Due to the local Harnack inequality, the behavior of the Green
function near the pole $x_0$, and \cite[Lemma 6.3]{P},  it follows
that there exists $C>0$ such that
$$C^{-1}G_{P_j}^\Omega(x,x_0)\leq u_j(x)\leq C G_{P_j}^\Omega(x,x_0),$$
for all $x\in \Gamma_\varepsilon$, and $j\geq 1$. Since $u_j$ and
$G_{P_j}^\Omega(x,x_0)$ are positive solutions of minimal growth
 of the operator $P_j$ in
a neighborhood of infinity in $\Omega\setminus \Gamma$ (see
\cite[Lemma 5.2]{P}),  it follows that \be\label{ujGj}
C^{-1}G_{P_j}^\Omega(x,x_0)\leq u_j(x)\leq CG_{P_j}^\Omega(x,x_0),
\ee for all $x\in \Omega\cap
\{\mbox{dist\,}(x,\Gamma)>\varepsilon\} $ and $j\geq 1$.

By taking a subsequence, we may assume that $\mu_j\to \mu_0$, and
$\{u_j\}$ converges in the open compact topology to a solution
$u\in \mathcal{C}_{\mu_0}(\Omega)$.  Clearly, $\mu_0\leq
\mu_\Omega$.

Since $\partial B_j$ are smooth, and $u_j$ vanish on $\partial
B_j\cap \partial \Omega_j$, it follows by \cite{D} and elliptic
regularity that $u$ vanishes on $\Gamma$.

Denote $P_0:=-\Delta-\mu_0 V$, and note that $P_0$ is subcritical
in $\Omega$. By (\ref{ujGj}) and the boundary Harnack principle,
$$C_1^{-1}G_{P_0}^\Omega(x,x_0)\leq u(x)\leq
C_1G_{P_0}^\Omega(x,x_0),$$ for all $x\in \Omega\setminus
\{\mbox{dist\,}(x,x_o)<\varepsilon/2\}$.
 Consequently, $u$ is a global positive
solution of the equation $P_0u=0$ in  $\Omega$ which has minimal
growth in a neighborhood of infinity in $\Omega$. In other words,
$u$ is a ground state of the operator $-\Delta-\mu_0 V$ in
$\Omega$. But this is a contradiction, since for $\mu\leq
\mu_\Omega$, the operator $-\Delta-\mu V$ is subcritical in
$\Omega$.
\end{proof}

We conclude this section with a direct proof of
Corollary~\ref{thm2old}.

\begin{proof}({\bf proof of Corollary~\ref{thm2old}})
Let $\mathrm{C}\varsubsetneq X$ be a Lipschitz cone, and let
$D\subset S_1^{N-1}$ be the Lipschitz domain so that
$\mathrm{C}=\{(r,\omega)\, | \,r\in(0,\infty), \omega\in D\}$. Let
$W\in L^p_{loc}(X)$ be a nonzero nonnegative function with a
compact support in $\mathrm{C}$, and set
$V(x)=\frac{1}{|x|^2}-W(x)$. Clearly, $\mu_\mathrm{C}=
\frac{(N-2)^2}{4}+\lambda_D$, where $\lambda_D$ is the Dirichlet
principal eigenvalue of $-\Delta_S$ on $D$. Moreover, the operator
$-\Delta-\mu_\mathrm{C}V$ is subcritical in $\mathrm{C}$. Let
$B_j\Subset X$ be a decreasing sequence of smooth domains, such
that $\mathrm{C}_j:=B_j\cup \mathrm{C}$ are connected for all
$j\geq 1$, and ${\rm
int}\left(\cap_j\mathrm{C}_j\right)=\mathrm{C}$. Fix $x_0\in
\mathrm{C}$.

Suppose that that $\mu_j:=\mu_{\mathrm{C}_j}$ is attained for all
$j\ge 1$, and let $u_j\in \mathcal{C}^0_{\mu_j}(\mathrm{C}_j) \cap
\merhav(\mathrm{C}_j) $ be the corresponding minimizer such that
$u_j(x_0)=1$. By Lemma \ref{lemminimizer}, $u_j$ is a positive
solution of the operator $-\Delta-\mu_j V$ of minimal growth in a
neighborhood of infinity in $\mathrm{C}_j$ .

We denote
$$\alpha_{j,\pm}:=\frac{-(N-2)\pm\sqrt{(N-2)^2-4(\mu_j-\lambda_D)}}{2},
\qquad  v_{j,\pm}(x):=r^{\alpha_{j,\pm}}v_D(\omega).$$ where $v_D$
is the Dirichlet principal eigenfunction of $-\Delta_S$ on $D$.
Fix $0<R_1<R_2$ such that $\supp W\subset \{R_1<|x|<R_2\}$ and
$$\mathrm{C}_j\cap
\left(\{|x|<R_1\}\cup\{|x|>R_2\}\right)\subset \mathrm{C}.$$ By
\cite[Theorem 6.3]{P}, there exists $C>0$ such that
$$C^{-1}v_{j,\pm}(x)\leq u_j(x)\leq Cv_{j,\pm}(x),$$
for all $x\in \mathrm{C}\cap
\left(\{|x|=R_1\}\cup\{|x|=R_2\}\right)$, and $j\geq 1$.

Since $u_j\in \merhav(\mathrm{C}_j)$, and $v_{j,+}$ is a positive
solution of minimal growth at the singular point $\zeta=0$
 of the operator $-\Delta-\frac{\mu_j}{|x|^2}$ in
$\mathrm{C}$, it follows that $u_j$ is a positive solution of
minimal growth at the singular point $\zeta=0$, and
\be\label{ujvjp} C^{-1}v_{j,+}(x)\leq u_j(x)\leq Cv_{j,+}(x), \ee
for all $x\in \mathrm{C}\cap \{|x|<R_1\}\}$ and $j\geq 1$.

Similarly, since $u_j\in \merhav(\mathrm{C}_j)$, and $v_{j,-}$ is
a positive solution of minimal growth at the singular point
$\zeta=\infty$ of the operator $-\Delta-\frac{\mu_j}{|x|^2}$ in
$\mathrm{C}$, it follows that $u_j$ is a positive solution of
minimal growth at the singular point $\zeta=\infty$, and
\be\label{ujvjm}C^{-1}v_{j,-}(x)\leq u_j(x)\leq Cv_{j,-}(x),\ee
for all $x\in \mathrm{C}\cap \{|x|>R_2\}\}$, and $j\geq 1$.

By taking a subsequence, we may assume that $\mu_j\to \mu_0$, and
$\{u_j\}$ converges in the open compact topology to a solution
$u\in \mathcal{C}_{\mu_0}(\mathrm{C})$. Moreover, since $B_j$ are
smooth, and $u_j$ vanish on $\partial B_j\cap
\partial \mathrm{C}_j$, it follows by \cite{D} and elliptic
regularity that $u$ vanishes on $\partial
\mathrm{C}\setminus\{0\}$.
 So, $u\in
\mathcal{C}^0_{\mu_0}(\mathrm{C})$.

Clearly, $\mu_0\leq \mu_\mathrm{C}$. Furthermore, $\alpha_{j,\pm}
\to \alpha_{0,\pm}:=
\frac{-(N-2)\pm\sqrt{(N-2)^2-4(\mu_0-\lambda_D)}}{2}$. Therefore,
$$C^{-1}|x|^{\alpha_{0,+}}v_D(\frac{x}{|x|})\leq u(x)\leq C|x|^{\alpha_{0,+}}v_D(\frac{x}{|x|}),$$
 for all $x\in
\mathrm{C}\cap \{|x|<R_1\}\}$, and
$$C^{-1}|x|^{\alpha_{0,-}}v_D(\frac{x}{|x|})\leq u(x)\leq C|x|^{\alpha_{0,-}}v_D(\frac{x}{|x|}),$$
 for all $x\in
\mathrm{C}\cap \{|x|>R_2\}\}$. Consequently, $u$ is a global
positive solution of the equation $(-\Delta-\mu_0 V)u=0$ in
$\mathrm{C}$ which has minimal growth in a neighborhood of
infinity in $\mathrm{C}$. In other words, $u$ is a ground state of
the operator $-\Delta-\mu_0 V$ in $\mathrm{C}$. But this is a
contradiction, since for  $\mu\leq \mu_\mathrm{C}$, the operator
$-\Delta-\mu V$ is subcritical in $\mathrm{C}$.
\end{proof}

\vskip 3mm

\begin{center}
{\bf Acknowledgments} \end{center} The authors wish to thank
M.~Murata and M.~Solomyak for valuable discussions. This research
was done at the Technion as K.~T. was a Lady Davis Visiting
Professor. K.~T. would like to thank M.~Marcus, Y.~Pinchover, and
G.~Wolansky for their hospitality. The work of Y.~P. was partially
supported by the RTN network ``Nonlinear Partial Differential
Equations Describing Front
 Propagation and Other Singular Phenomena", HPRN-CT-2002-00274,
and the Fund for the Promotion of Research at the Technion. The
work of K.~T. was partially supported by the Swedish Research
Council.

\end{document}